\documentclass[12pt]{article}

\usepackage{amsfonts}
\usepackage{amsmath}
\usepackage{amssymb}
\usepackage{indentfirst}

\usepackage[english,russian]{babel}

\renewcommand{\le}{\leqslant}

\begin{document}
\title{On the coloring of 3-element subsets\footnote{This paper is prepared under the supervision of A.M.Raigorodskiy and is submitted to the Moscow Mathematical Conference for High-School Students. Readers are invited to send their remarks and reports on this paper to mmks@mccme.ru}}
\author{D. Zakharov}
\date{}
\maketitle
\newtheorem{theorem}{Theorem}
\newtheorem{prop}{Proposition}
\newtheorem{rem}{Remark}

\renewcommand{\refname}{References}
\begin{theorem}\label{t1}
Let $p=8k-1$ be a prime. Then we can color the set ${R_{p+2} \choose 3}$ of 3-element subsets of $R_{p+2}=\{1, \ldots, p+2\}$ into $p$ colors such that any two sets $x, y \in {R_{p+2} \choose 3}, ~ |x\cap y| = 2$ have different colors.
\end{theorem}

Professor Raigorodskiy confirmed that this result is new.

\textbf{Remarks.}  (a) Obviously, the set $R_{p+2} \choose 3$ can not be colored in a fewer than $p$ colors because in the set 
$$\{ \{1, 2, 3\}, \{1,2,4\}, \ldots, \{1, 2, p+2\} \}$$
any pair of triplets must have differrent colors.

(b) It is clear from the proof that the theorem also holds for primes $p$ such that  $-1 \not\equiv 2^r \pmod p$ for each $r$. 

(c) We can consider graph $G(n, 3, 2)$ (see [1-6]) which has vertices set $R_n \choose 3$ that are connected if their intersection has two elements, i.e.
$$
G(n, 3, 2) = (V, E),~ V={R_n \choose 3}, ~E=\{\{x, y\}: |x\cap y|= 2\}
$$
In this terminology the theorem states that $\chi(G(p+2, 3, 2))\le p$ (and so $\chi(G(p+2, 3, 2))=p$).

\textbf{Aknowledgements.} I would like to acknowledge my advisor prof. A. M.Raigorodskiy for his attention to this work.

\section{The proof of theorem~\ref{t1}}

\subsection{Construction of coloring}

\newtheorem{lemma}{Lemma}

Let $V = {R_n \choose 3}, N=\{p+1, p+2\}$. Consider 4 sets:
$$V_0=\{v\in V, v \cap N = \emptyset\},$$
$$V_2=\{v\in V, N \subset v\},$$
$$W_1=\{v\in V, v \cap N = \{p+1\} \},$$
$$W_2=\{v\in V, v \cap N = \{p+2\} \}$$

\noindent Every vertex  $x=\{x_1, x_2, x_3\} \in V_0$ we paint in color 
$$
c(x)=x_1+x_2+x_3 \pmod p,
$$
 vertex $x=(x_1, p+2, p+1)\in V_2$ we paint in color 
 $$
 c(x)=3x_1 \pmod p.
$$
For the coloring of $W_1$ and $W_2$ we will use the following lemma, which will be proved in 2.3.

Let $\widetilde R_p=\{(x, y)|x, y \in R_p, x \not = y\}$.

\begin{lemma}\label{l}
There exists a map $f: \widetilde{R_p} \rightarrow \mathbb{Z}_2$ such that 
\begin{enumerate}

\item $f(x,y)\ne f(y,x)$ for each $x,y$ 
\item $f(x,y)\ne f(\frac{x+y}{2},x)$ for each $x,y$
(here the division is the division in $\mathbb{Z}_p$), 
\end{enumerate}

\end{lemma}

Define a function $f_1:\widetilde{R_p}\to\mathbb{Z}_p$ by $f_1(x, y) = x$ if $f(x,y)=0$ and $f_1(x, y) = y$ if $f(x,y)=1$.  
Define a function $f_2:\widetilde{R_p}\to\mathbb{Z}_p$ by $f_2(x, y) = x+y-f_1(x,y)$.
Now let us paint vertex 
$$
x=(x_1, x_2, p+i)\in W_i
$$
in color
$$
c(x)=x_1+x_2+f_i(x_1, x_2) \pmod p.
$$ 

Obviuosly, we have constructed a coloring of all elements of $R_{p+2} \choose 3$ in $ p $ colors. 

\subsection{The proof that coloring is regular}

Let us take any two elements 
$$
x = (x_1, x_2, x_3), ~~~ y = (x_1, x_2, x_4) \in {R_{p+2} \choose 3}
$$
and consider following cases:

\paragraph{Case 1: $x, y \in V_0$.} Then, obviously, 
$$
c(x) \equiv x_1+x_2+x_3 \not \equiv x_1+x_2+x_4 \equiv c(y) \pmod p.
$$

\paragraph{Case 2: $x, y \in V_2,~x_1, x_2 \in N$}. Then 
$$
c(x) \equiv 3x_3\not \equiv 3x_4 \equiv c(y) \pmod p,
$$
because $ p > 3$ is prime, i.e. $p$ is not divisible by 3.

\paragraph{Case 3: $x \in V_0, $ $ y \in V_2$.} This case is impossible.

\paragraph{Case 4: $x \in W_i, $ $ y \in V_2$.} We can assume that $x_1 = p+2 $, $ x_4 = p+1$, so
$$
c(x)-c(y)\equiv x_2+x_3+f_i(x_2, x_3)-3x_2 \equiv x_3 + f_i(x_2, x_3) - 2x_2 \pmod p.
$$ 
If $f_i(x_2, x_3)=x_2$, then
$$
c(x)-c(y) \equiv x_3 - x_2 \not \equiv 0 \pmod p,
$$ 
else $f_i(x_2, x_3)=x_3$ and
$$c(x)-c(y) \equiv 2x_3 - 2x_2 \not \equiv 0 \pmod p,$$

\paragraph{Case 5: $x\in W_i $, $y \in V_0, x_3 \in N$.} We have 
$$
c(x)-c(y) \equiv x_1+x_2+f_i(x_1, x_2)-x_1-x_2-x_4 \equiv f_i(x_1, x_2)-x_4 \not \equiv 0 \pmod p
$$
by definition of $f_i$.

\paragraph{Case 6: $x \in W_1$, $y \in W_2$.} In this case $x_3=p+1$, $x_4=p+2$, consequently
$$
c(x)-c(y)\equiv f_1(x_1, x_2)-f_2(x_1, x_2) \not \equiv 0 \pmod p
$$
by definition of $f_1$ and $f_2$. 

\paragraph{Case 7: $x, y \in W_1, x_1 = p+1$.} We can write that
$$
c(x)-c(y) \equiv x_2+x_3+f_1(x_2, x_3)-x_2-x_4-f_1(x_2, x_4)\equiv x_3+f_1(x_2, x_3)-x_4-f_1(x_2, x_4) \pmod p. 
$$
Consider subcases.

\paragraph{Subcase 7.1: $f_1(x_2, x_3)=f_1(x_2, x_4)=x_2$.} Then
$$
x_3+f_1(x_2, x_3)-x_4-f_1(x_2, x_4)\equiv x_3-x_4 \not \equiv 0 \pmod p.
$$

\paragraph{Subcase 7.2: $f_1(x_2, x_3)=x_3, $ $ f_1(x_2, x_4)=x_4$.} Then
$$
x_3+f_1(x_2, x_3)-x_4-f_1(x_2, x_4)\equiv 2x_3-2x_4 \not \equiv 0 \pmod p.
$$

\paragraph{Subcase 7.3: $f_1(x_2, x_3)=x_2, $ $ f_1(x_2, x_4)=x_4$.} Then
$$
x_3+f_1(x_2, x_3)-x_4-f_1(x_2, x_4)\equiv x_3+x_2-2x_4 \pmod p.
$$
Suppose $x_3+x_2\equiv 2x_4 \pmod p$. So we have $x_4=\frac{x_2+x_3}{2}$. From this we get 
$$
f_1(x_2, x_3)=x_2, ~~~ f_1\left(x_2, \frac{x_2+x_3}{2}\right)=\frac{x_2+x_3}{2}.
$$
And so 
$$
f(x_2, x_3)=1, ~f\left(x_2, \frac{x_2+x_3}{2}\right)=0 \Rightarrow f(x_3, x_2) = 0 
$$
which contradicts to the properties of $f$. Hence, $c(x)-c(y)\not\equiv 0 \pmod p$.

\vskip+0.3cm

Thus, we considered all cases and the constructed coloring is regular.

\subsection{The proof of lemma \ref{l}}

\begin{lemma} \label{l2}
Let $p=8k-1$ be a prime. Then $-1 \not\equiv 2^r \pmod p$ for each $r$. 
\end{lemma}

\textbf{Proof.} Denote by $d$ the order of 2 modulo $p$. We denote by $\left(\frac{a}{p}\right)$ the Legendre symbol. It's known, that 
$$
\left(\frac{2}{p}\right)=(-1)^{\frac{p^2-1}{8}}=(-1)^{\frac{64k^2-16k+1-1}{8}}=1.
$$
Thus, 2 is a quadratic residue in $\mathbb{Z}_p$. Then $2^{\frac{p-1}{2}}\equiv 1 \pmod p$ and so $d|\frac{p-1}{2}=4k-1$, i.e. $d$ is odd. 

Suppose, that there exists minimal number $l$, such that $2^l\equiv -1 \pmod p$.  
Then $d | 2l$ but $d$ is odd, whence $d|l$ and $2^l \equiv 1 \pmod p$. Lemma \ref{l2} is proved.

\vskip+1cm

Let us define a graph $G$ with vertex set $\widetilde{R_p}$ and in which we connect $(x, y)$ with $$
(y, x),~ (\frac{x+y}{2}, x), ~ (y, 2x-y)
$$

It is easy to see that if $(x, y)$ and $(a, b)$ are connected then $\frac{x-y}{a-b}=-2^s$ for some $s$. So if there is an odd cycle $(x_1, y_1), \ldots, (x_l, y_l)$ then 
$$
x_1-y_1 \equiv -2^{s_2} (x_2-y_2) \equiv \ldots \equiv (-1)^{l-1}2^{s_2+\ldots+s_l}(x_n-y_n)\equiv (-1)^l 2^{s_1+\ldots+s_l}(x_1-y_1)
$$
And we get
$$
(x_1-y_1)(2^S+1)\equiv 0 \pmod p
$$
which contradicts to Lemma \ref{l2}. Then $G$ has not odd cycles so it is bipartite. 
\vskip+1cm

Now let us take some 2-coloring of $G$ ~$V=\mathcal M_1 \cup \mathcal M_2$. We define the map $f$ as:
$$
f(x, y) = 1  \Leftrightarrow (x, y) \in \mathcal M_1
$$
Clearly $f(x, y) \not = f(y, x)$ and $f(x, y) \not = f(\frac{x+y}{2}, x)$ because corresponding vertices are adjacent. Lemma \ref{l} is proved.


\begin{thebibliography} {20}

\bibitem{Rai1} A.M. Raigorodskii, {\it Cliques and cycles in distance graphs and graphs of diameters}, ``Discrete Geometry and Algebraic Combinatorics'',
AMS, Contemporary Mathematics, 625 (2014), 93 - 109.

\bibitem{Rai2} A.M. Raigorodskii, {\it Coloring Distance Graphs and Graphs of Diameters}, Thirty Essays on Geometric Graph Theory, J. Pach ed., Springer, 2013, 
429 - 460.

\bibitem{Rai3} B. Bollob\'as, B.P. Narayanan, A.M. Raigorodskii, {\it On the stability of the Erd\H{o}s--Ko--Rado theorem}, J. Comb. Th. Ser. A, 137 (2016), 64 - 78.

\bibitem{Rai5} A.M. Raigorodskii, {\it Combinatorial geometry and coding theory}, Fundamenta Informatica, 145 (2016), 359 - 369.

\bibitem{BKK} A.V. Bobu, O.A. Kostina, A.E. Kupriyanov, {\it Independence numbers and chromatic numbers of some distance graphs}, Problemy Peredachi Informatsii, 2015, Vol. 51, No. 2, pp. 86–98.

\bibitem{Rai8} J. Balogh, A.V. Kostochka, A.M. Raigorodskii, {\it Coloring some finite sets in $ {\mathbb R}^n $}, Discussiones Mathematicae
Graph Theory, 33 (2013), N1, 25 - 31.

\end{thebibliography}
\end{document}